\documentclass[a4paper,11pt]{article}
\usepackage{amssymb}
\usepackage{amsfonts}
\usepackage{amsmath}

\input{tcilatex}

\begin{document}

\title{Circle actions and Suspension operations on Smooth manifolds}
\author{Haibao Duan \thanks{%
Supported by National Science foundation of China No. 12331003}}
\date{}
\maketitle

\begin{abstract}
Let $M$ be a smooth manifold with $\dim M\geq 3$ and a base point $x_{0}$.
Surgeries along the oriented circle $S^{1}\times \{x_{0}\}$ on the product $%
S^{1}\times M$ yield two manifolds $\Sigma _{0}M$ and $\Sigma _{1}M$, called
the suspensions of $M$.

The suspension operations play a basic role in the construction and
classification of the smooth manifolds which admit free circle actions. This
paper is devoted to present such results and evidences.
\end{abstract}

\section{Introduction and main results}

Let $\mathcal{T}_{n}$ be the set of diffeomorphism classes of closed,
oriented and smooth manifolds of dimension $n$. A free circle action $%
S^{1}\times M\rightarrow M$ on a manifold $M\in \mathcal{T}_{n}$ is called 
\textsl{regular} if the quotient space $M/S^{1}$, equipped with the induced
topology and orientation, belongs to $\mathcal{T}_{n-1}$. For a manifold $%
M\in \mathcal{T}_{n}$ we would like to known whether there is a regular
circle action on $M$. If the answer is affirmative, we would like to decide
furthermore the diffeomorphism types of the quotients $M/S^{1}$. If $M$ is a
homotopy sphere, these problems were raised and studied by Montgemory-Yang,
Edmonds, Schultz, Madsen-Milgram and Wall \cite{MY,E,Sc,MM,Wall}. In the
case $M$ is an $1$-connected $5$-manifolds, the problems was solved by Duan
and Liang in \cite{DL}.

This paper introduces a general method to construct smooth manifolds $M$
that admit regular circle actions, while the diffeomorphism types of the
quotients $M/S^{1}$ can be made explicit. The central ingredient of this
construction is the \textsl{suspension operations }$\Sigma _{i}:\mathcal{T}%
_{n}\rightarrow \mathcal{T}_{n+1}$ (with $i=0,1$) we are about to describe.
In this paper $D^{n}$ denotes the unit disc on the Euclidean $n$-space $%
\mathbb{R}^{n}$, with boundary the $(n-1)$-sphere $S^{n-1}=\partial D^{n}$.
We shall assume that $n\geq 3$, unless otherwise stated.

The fundamental group $\pi _{1}(SO(n))$ of the special orthogonal group $%
SO(n)$ with order $n\geq 3$ is isomorphic to $\mathbb{Z}_{2}$ with the
generator

\begin{quote}
$\alpha (t)=\left( 
\begin{array}{cc}
\cos 2\pi t & -\sin 2\pi t \\ 
\sin 2\pi t & \cos 2\pi t%
\end{array}%
\right) \oplus I_{n-2}$, $t\in \lbrack 0,1]$,
\end{quote}

\noindent where $I_{r}$ denotes the identity matrix of rank $r$. It gives
rise to the diffeomorphism of the product $S^{1}\times D^{n}$

\begin{enumerate}
\item[(1.1)] $\widetilde{\alpha }:S^{1}\times D^{n}\rightarrow S^{1}\times
D^{n}$ by $\widetilde{\alpha }(t,x)=(t,\alpha (t)x)$,
\end{enumerate}

\noindent whose restriction on the boundary $S^{1}\times S^{n-1}$ is denoted
by $\tau $. Recall that any orientation preserving diffeomorphism of $%
S^{1}\times S^{n-1}$ is isotopic either to identity $id$ or to $\tau $ (e.g. 
\cite[p.241]{Wall}). For a manifold $N\in \mathcal{T}_{n}$ with a base point 
$x_{0}\in N$, and an orientation preserving smooth embedding $%
(D^{n},0)\subset (N,x_{0})$ centered at $x_{0}$, we introduce two adjoint
manifolds $\Sigma _{i}N$, $i=0,1$, by

\begin{enumerate}
\item[(1.2)] $\Sigma _{0}N:=D^{2}\times S^{n-1}\cup _{id}(S^{1}\times (N-%
\overset{\circ }{D^{n}}))$,

$\Sigma _{1}N:=D^{2}\times S^{n-1}\cup _{\tau }(S^{1}\times (N-\overset{%
\circ }{D^{n}}))$,
\end{enumerate}

\noindent respectively, where $\cup _{f}$ means gluing the boundaries
together by the diffeomorphism $f$. Alternatively, the manifolds $\Sigma
_{0}N$ and $\Sigma _{1}N$ are obtained by surgeries along the oriented
circle $S^{1}\times \{x_{0}\}$ on $S^{1}\times N$ with the framings $id$ and 
$\widetilde{\alpha }$, respectively. By the homogeneity lemma \cite[\S 4]%
{Mil1} the diffeomorphism types of $\Sigma _{i}N$ are irrelevant to the
choice of a smooth embedding $(D^{n},0)\subset (N,x_{0})$. This justifies
the following notion.

\bigskip

\noindent \textbf{Definition 1.1. }The correspondences $\Sigma _{i}:\mathcal{%
T}_{n}\rightarrow \mathcal{T}_{n+1}$ by $N\rightarrow $ $\Sigma _{i}N$ ($%
i=0,1$) are called \textsl{the suspension operations on }$\mathcal{T}_{n}$%
.\hfill $\square $

\bigskip

The relationship between the two suspensions $\Sigma _{0}N$ and $\Sigma
_{1}N $ is quite subtle. We can only derive partial information about it.
Let $Diff(N,x_{0})$ be the group of diffeomorphisms of $N\in \mathcal{T}_{n}$
that fix the base point $x_{0}$, and let $GL(n)$ be the general linear group
of order $n$. Sending a diffeomorphism $h\in Diff(N,x_{0})$ to the tangent
map of $h$ at $x_{0}$ yields a homomorphism

\begin{quote}
$e_{N}:Diff(N,x_{0})\rightarrow GL(n)$, $h\rightarrow T_{x_{0}}h$,
\end{quote}

\noindent called \textsl{the linear representation} of the group $%
Diff(N,x_{0})$. A manifold $N\in \mathcal{T}_{n}$ ($n\geq 3$) is called $%
\Sigma $\textsl{-stable} if the induced action on the fundamental groups

\begin{quote}
$e_{N\ast }:\pi _{1}(Diff(N,x_{0}),id)\rightarrow \pi _{1}(GL(n),I_{n})=%
\mathbb{Z}_{2}$
\end{quote}

\noindent surjects. The following result will be shown in Section 3.

\bigskip

\noindent \textbf{Theorem A.} \textsl{If }$N\in \mathcal{T}_{n}$ \textsl{is} 
$\Sigma $\textsl{-stable, then }$\Sigma _{0}N$ $=\Sigma _{1}N$\textsl{.}

\bigskip

As supplements to Theorem A we shall see that the product $S^{k}\times M$
with $k\geq 3$ are all $\Sigma $-stable. On the contrary, $\Sigma _{1}M$%
\textsl{\ }$\neq \Sigma _{0}M$ for the quaternionic projective plane $M=%
\mathbb{H}P^{2}$.

For a regular circle action on a manifold $E\in \mathcal{T}_{n}$ we can
regard the quotient map

\begin{quote}
$p:E\rightarrow B:=E/S^{1}$
\end{quote}

\noindent as an oriented circle bundle over $B$. For a manifold $N\in 
\mathcal{T}_{n-1}$ let $B\#N$ be the connected sum of $B$ and $N$, let $%
f_{N}:B\#N\rightarrow B$ be the canonical projection. Denote the induced
bundle $f_{N}^{\ast }E$ by

\begin{enumerate}
\item[(1.3)] $p_{N}:E_{N}\rightarrow $ $B\#N$.
\end{enumerate}

\noindent Our main result expresses the total space $E_{N}$ in term of $E$
and $N$, where $w_{2}(M)$ stands for the\ second Stiefel--Whitney class\ of
a manifold $M$.

\bigskip

\noindent \textbf{Theorem B.} \textsl{Let }$S^{1}\times E\rightarrow E$ 
\textsl{be a regular circle action on an }$1$\textsl{-connected manifold }$%
E\in \mathcal{T}_{n}$\textsl{\ with\ quotient map }$p:E\rightarrow B$\textsl{%
, where} $n\geq 5$\textsl{. }

\textsl{Then for any }$N\in \mathcal{T}_{n-1}$ \textsl{we have}

\begin{enumerate}
\item[(1.4)] $E_{N}=\left\{ 
\begin{tabular}{l}
$E\#\Sigma _{0}N$ \textsl{if} $w_{2}(B)\neq 0$\textsl{;} \\ 
$E\#\Sigma _{1}N$ \textsl{if} $w_{2}(B)=0$\textsl{.}%
\end{tabular}%
\right. $
\end{enumerate}

\noindent \textsl{In addition, if }$p^{\ast }(w_{2}(B))=w_{2}(E)\neq 0$%
\textsl{,\ then}

\begin{enumerate}
\item[(1.5)] $E_{N}=E\#\Sigma _{0}N=E\#\Sigma _{1}N$\textsl{.}
\end{enumerate}

In Theorem B the $1$-connectness restriction on $E$ can not be excluded. For
an integer $q\geq 2$ let $p:E\rightarrow \mathbb{C}P^{2}$ be the oriented
circle bundle over the complex projective plane with Euler class $q\cdot x$,
where $x\in H^{2}(\mathbb{C}P^{2})$ is a generator. Then $E$ is the $5$%
-dimensional lens space $L^{5}(q)$ whose fundamental group is cyclic with
order $q$. It can be shown that $E_{N}$ is homeomorphic neither to $%
E\#\Sigma _{0}N$ nor to $E\#\Sigma _{1}N$, unless $N\in \mathcal{T}_{4}$ is
the sphere $S^{4}$.

Theorem B is useful to produce manifolds that admit regular circle actions
with explicit quotients. We present such examples.

\bigskip

\noindent \textbf{Corollary 1.3. }\textsl{Let}\textbf{\ }$E\in \mathcal{T}%
_{n}$ \textsl{be an }$1$\textsl{-connected manifold that admits a regular
circle action with\ quotient }$B$\textsl{, where} $n\geq 5$\textsl{.}

\textsl{Then, for any }$N\in \mathcal{T}_{n-1}$\textsl{\ one of the two
connected sums }$E\#\Sigma _{i}N$\textsl{,} $i=0,1$\textsl{,\ admits a
regular circle action with quotient }$B\#N$\textsl{.}\hfill $\square $

\bigskip

\noindent \textbf{Example 1.4. }Let $p:S^{2n+1}\rightarrow \mathbb{C}P^{n}$
be the Hopf fibration over the complex projective $n$-space $\mathbb{C}P^{n}$%
, $n\geq 2$. By Corollary 1.3, for any $N\in \mathcal{T}_{2n}$ either $%
\Sigma _{0}N$ or $\Sigma _{1}N$ admits a regular circle action with\
quotient $\mathbb{C}P^{n}\#N$. In particular, by the homotopy exact sequence
of circle fibrations one concludes that the homotopy groups of the
suspension $\Sigma _{i}N$\ are

\begin{enumerate}
\item[(1.6)] $\pi _{r}(\Sigma _{i}N)=\left\{ 
\begin{array}{c}
\pi _{r}(N)\text{ }\QTR{sl}{if}\text{ }r=1\text{ }\QTR{sl}{or}\text{ }2%
\QTR{sl}{,} \\ 
\pi _{r}(\mathbb{C}P^{n}\#N)\text{ }\QTR{sl}{if}\text{ }r\geq 3\QTR{sl}{,}%
\end{array}%
\right. $
\end{enumerate}

\noindent where $n\equiv i\func{mod}2$.

In general, by a \textsl{fake complex projective }$n$\textsl{-space} we mean
a smooth manifold $F_{\mathbb{C}P^{n}}$ that is homotopy equivalent to the
complex projective space $\mathbb{C}P^{n}$, $n\geq 2$. For the
classifications of such manifolds, see \cite[Theorem 4.9]{MM} and \cite[%
Chapter 14C]{Wall}. In particular, if we let $p:E\rightarrow F_{\mathbb{C}%
P^{n}}$ be the oriented circle bundle on $F_{\mathbb{C}P^{n}}$ whose Euler
class generates the group $H^{2}(F_{\mathbb{C}P^{n}})=\mathbb{Z}$, then $E$
is a homotopy $(2n+1)$-spheres. In addition, since the second
Stiefel-Whitney class is a homotopy invariant of smooth manifolds, we have $%
w_{2}(F_{\mathbb{C}P^{n}})\neq 0$ if and only if $n$ even. Thus, Theorem B
implies that, for any fake complex projective space $F_{\mathbb{C}P^{n}}$
and any $N\in \mathcal{T}_{2n}$, the (topological) manifold $\Sigma _{i}N$
admits a smooth structure, together with a regular circle action, whose
quotient is diffeomorphic to $F_{\mathbb{C}P^{n}}\#N$, where $n\equiv i\func{%
mod}2$.\hfill $\square $

\bigskip

\noindent \textbf{Example 1.5.} For an integer $n\geq 2$ let $S^{2}%
\widetilde{\times }S^{n}$ be the only nontrivial smooth $S^{n}$-bundle over $%
S^{2}$. Take $B:=S^{2}\times S^{n-2}$ or $S^{2}\widetilde{\times }S^{n-2}$
with $n\geq 4$, and let $p:E\rightarrow B$ be the circle bundle over $B$
whose Euler class is the generator of $H^{2}(B)$ corresponding to the
orientation class of the factor $S^{2}\subset B$. Since the total space $%
E=S^{3}\times S^{n-2}$ is $1$-connected, $w_{2}(E)=0$, and since $%
w_{2}(B)\neq 0$ if and only if $B=S^{2}\widetilde{\times }S^{n-2}$, we
conclude by Corollary 1.3 that for any $N\in \mathcal{T}_{n-1}$\ with $n\geq
5$, both manifolds $S^{3}\times S^{n-2}\#\Sigma _{i}N$ admit regular circle
actions with quotients $S^{2}\widetilde{\times }S^{n-2}\#N$ and $S^{2}\times
S^{n-2}\#N$, respectively.\hfill $\square $

\bigskip

Combined with general properties of the suspensions $\Sigma _{i}$ to be
developed in Section 4, Theorem B is applicable to study the classification
problem of the manifolds that admit regular circle actions. In the 1960's
Smale and Barden \cite{B,S} classified the 1-connected 5-manifolds in term
of a collection of five basic ones $\{S^{2}\times S^{3},S^{2}\widetilde{%
\times }S^{3},W,M_{k},X_{2^{i}}\}$, where $k\geq 2,i\geq 1$, see Section 5
for detailed constructions of these manifolds. On the other hand, in \cite%
{GL1,L,GL} R. Goldstein and L. Lininger initiated the problem of classifying
the 1-connected 6-manifolds that admit regular circle actions. Our next
result completes the project (compare with \cite[Theorem 4, Theorem 6]{L}).

\bigskip

\noindent \textbf{Theorem C. }\textsl{If }$M$\textsl{\ is an} $1$\textsl{%
-connected }$6$\textsl{-manifold that admits regular circle actions, then}

\begin{enumerate}
\item[(1.7)] $M=\left\{ 
\begin{tabular}{l}
$S^{3}\times S^{3}\#_{r}\Sigma _{0}(S^{2}\times S^{3})\#_{1\leq j\leq
t}\Sigma _{1}M_{k_{j}}\#\Sigma _{1}H$ \textsl{if} $w_{2}(M)\neq 0$\textsl{,}
\\ 
$S^{3}\times S^{3}\#_{r}\Sigma _{0}(S^{2}\times S^{3})\#_{1\leq j\leq
t}\Sigma _{i}M_{k_{j}}$ \textsl{if} $w_{2}(M)=0,$%
\end{tabular}%
\right. $
\end{enumerate}

\noindent \textsl{where }

\begin{quote}
$H\in \{S^{2}\widetilde{\times }S^{3},W,X_{k}\}$\textsl{, }$\Sigma
_{0}(S^{2}\times S^{3})=S^{3}\times S^{3}\#S^{2}\times S^{4}$\textsl{, }
\end{quote}

\noindent \textsl{and where the notation }$\#_{r}N$ \textsl{means} \textsl{%
the connected sum of }$r$\textsl{-copies of} $N$\textsl{\ (and so forth).}

\bigskip

Let $T^{k}=S^{1}\times \cdots \times S^{1}$ ($k$-copies) be the $k$%
-dimensional torus group with classifying space $BT^{k}$. For a
(topological) manifold $M$, and a set of cohomology classes $\left\{ \alpha
_{1},\cdots ,\alpha _{k}\right\} \subset H^{2}(M)$ of degree $2$, let $%
f_{i}:M\rightarrow \mathbb{C}P^{\infty }$ be the classifying map of $\alpha
_{i}$, and let

\begin{enumerate}
\item[(1.8)] $\pi :M(\alpha _{1},\cdots ,\alpha _{k})\rightarrow M$
\end{enumerate}

\noindent be the principal $T^{k}$-bundle on $M$ whose classifying map is
the product

\begin{quote}
$f=\left( f_{1},\cdots ,f_{k}\right) :M\rightarrow BT^{k}=\mathbb{C}%
P^{\infty }\times \cdots \times \mathbb{C}P^{\infty }$ ($k$-copies).
\end{quote}

\noindent Theorem B can be iterated to calculate the homeomorphism type of
the total space $M(\alpha _{1},\cdots ,\alpha _{k})$. We illustrate the
calculation for the cases where $M$ is an $1$-connected $4$-manifold, which
is not necessarily smoothable.

For an integer $k\geq 1$ define the sequence of integers $(b_{1},\cdots ,b_{%
\left[ \frac{k}{2}\right] })$ by

\begin{quote}
$b_{i}:=(k-1)\binom{k-1}{i}-\binom{k-1}{i+1}+(k-1)\binom{k-1}{i-1}-\binom{k-1%
}{i-2}$, $1\leq i\leq \left[ \frac{k}{2}\right] $.
\end{quote}

\noindent Accordingly, let $Q_{k}$ be the connected sum of the products of
spheres

\begin{enumerate}
\item[(1.9)] $Q_{k}:=\#_{c_{1}}(S^{{\small 3}}\times S^{{\small k+1}%
})\#_{c_{2}}(S^{{\small 4}}\times S^{{\small k}})\#_{c_{3}}\cdots
\#_{c_{r}}(S^{r{\small +2}}\times S^{{\small k-}r{\small +2}})$, $r=\left[ 
\frac{k}{2}\right] $,
\end{enumerate}

\noindent where $c_{i}=b_{i\text{ }}$with the only exception that $c_{r}=%
\frac{1}{2}b_{r}$ if $k\equiv 0\func{mod}2$.

Orlik and Raymond \cite[p.553]{OR} (resp. Davis and Januszkiewicz \cite[%
Example 1.20]{DJ})\textbf{\ }proved that\textbf{\ }any \textsl{toric }$4$%
\textsl{-manifold} is homeomorphic to a connected sum of copies of $\mathbb{C%
}P^{2},\overline{\mathbb{C}P}^{2}$ and $S^{2}\times S^{2}$. Buchstaber and
Panov \cite[Theorem 4.6.12]{BP} have further shown that (e.g. \cite[Theorem
6.3]{BM}, \cite{M}), if $M$ is a toric $4$-manifold and if $\left\{ \alpha
_{1},\cdots ,\alpha _{k}\right\} $ is the canonical basis of $H^{2}(M)$,
then the manifold $M(\alpha _{{\small 1}},\cdots ,\alpha _{{\small k}})$ is
diffeomorphic\textsl{\ }to\textsl{\ }$Q_{k}$. Indeed, this phenomenon holds
for all $1$-connected topological $4$-manifolds. Let $\beta _{2}(X)$ denote
the second Betti number of a CW-complex $X$.

\bigskip

\noindent \textbf{Theorem D. }\textsl{If }$M$\textsl{\ is an }$1$\textsl{%
-connected }$4$\textsl{-manifold} \textsl{and if }$\left\{ \alpha
_{1},\cdots ,\alpha _{k}\right\} $\textsl{\ is a basis of }$H^{2}(M)$\textsl{%
,} $k=\beta _{2}(X)$\textsl{, then }$M(\alpha _{{\small 1}},\cdots ,\alpha _{%
{\small k}})$ \textsl{is diffeomorphic to }$Q_{k}$.

\textsl{In particular,} \textsl{any }$1$\textsl{-connected }$4$\textsl{%
-manifold }$M$\textsl{\ with }$k=\beta _{2}(M)$ \textsl{is the quotient of a
(topological) principal }$T^{k}$\textsl{-action on} $Q_{k}$\textsl{.}

\bigskip

From the homotopy exact sequence of the bundle $\pi $ one finds that

\begin{quote}
$\pi _{r}(M(\alpha _{{\small 1}},\cdots ,\alpha _{{\small k}}))=\left\{ 
\begin{tabular}{l}
$0$ for $r=0,1,2$; \\ 
$\pi _{r}(M)$ for $r\geq 3$ (induced by $\pi _{\ast }$).%
\end{tabular}%
\right. $
\end{quote}

\noindent From $M(\alpha _{{\small 1}},\cdots ,\alpha _{{\small k}})\cong
Q_{k}$ by Theorem C we get

\bigskip

\noindent \textbf{Corollary 1.6.} \textsl{The }$2$\textsl{--connected cover
of any }$1$\textsl{-connected }$4$\textsl{-manifold }$M$ \textsl{is
homeomorphic to} $Q_{k}$\textsl{,} \textsl{where} $k=\beta _{2}(M)$\textsl{.}

\textsl{In particular, for }$3\leq r\leq k+3$\textsl{\ the homotopy group }$%
\pi _{r}(M)$\textsl{\ can be expressed in terms of that of spheres by the
formula}

\begin{enumerate}
\item[(1.10)] $\pi _{r}(M)=\pi _{r}(\underset{{\small c}_{{\small 1}}}{%
{\small \vee }}(S^{{\small 3}}\vee S^{{\small k+1}})\underset{{\small c}_{2}}%
{{\small \vee }}(S^{{\small 4}}\vee S^{{\small k}})\underset{{\small c}_{%
{\small 3}}}{{\small \vee }}\cdots \underset{{\small c}_{r}}{{\small \vee }}%
(S^{r{\small +2}}\vee S^{{\small k-}r{\small +2}}))$\textsl{,}
\end{enumerate}

\noindent \textsl{as well as} \textsl{the Hilton-Milnor's formula} \cite[%
p.533]{Wh}\textsl{.}\hfill $\square $

\bigskip

The map $\pi $ in (1.8) induces the fibration $\Omega \pi :\Omega M(\alpha _{%
{\small 1}},\cdots ,\alpha _{{\small k}})\rightarrow \Omega M$ in the loop
spaces, where $\Omega M(\alpha _{{\small 1}},\cdots ,\alpha _{{\small k}%
})\cong \Omega Q_{k}$ by Theorem D. Since the space $\Omega Q_{k}$ is $1$%
-connected, while the fiber $\Omega T^{k}$ is homotopy equivalent to the
Eilenberg-McLane space $K(\oplus _{k}\mathbb{Z},0)$, we get

\bigskip

\noindent \textbf{Corollary 1.7.}\textsl{\ For any} $1$\textsl{-connected }$4
$\textsl{-manifold }$M$ \textsl{the universal cover of the loop space} $%
\Omega M$ \textsl{is }$\Omega Q_{k}$\textsl{, where }$k=\beta _{2}(M)$%
\textsl{.}\hfill $\square $

\bigskip

The remaining sections of the paper are organized as follows. Section 2
introduces isotopy invariants for the framed circles on a manifold, and
establishes Theorem B. Section 3 is devoted to show Theorem A, and present
relevant examples. Based on general properties of the suspension operators $%
\Sigma _{i}$ developed in Section 4, the proofs of Theorem C and D are
completed in Sections 5 and 6, respectively.

\bigskip

\noindent \textbf{Acknowledgements.} This paper is based on a seminar talk
given at Institute of Mathematics, CAS. The author would like to thank Dr.
Yang Su, Feifei Fan, Xueqi Wang for useful communication. Thanks also due to
my referee for pointing out an error in an earlier version of the paper.

\section{Invariants of framed circles on a manifold.}

For an orientation preserving smooth embedding $F:S^{1}\times
D^{n-1}\rightarrow M$ into a manifold $M\in \mathcal{T}_{n}$ we set $%
S^{1}=F(S^{1}\times \{0\})\subset M$, and call the pair $(S^{1},F)$ a 
\textsl{framed} \textsl{circle} on $M$. Two framed circles $(S^{1},F)$ and $%
(S^{1},F^{\prime })$ are called \textsl{isotopic} if there is a smooth
embedding $H:(S^{1}\times D^{n-1})\times I\rightarrow M\times I$ such that $%
H((S^{1}\times D^{n-1})\times \{t\})\subset M\times \{t\}$ and

\begin{quote}
$H(x,0)=(F(x),0)$, $H(x,1)=(F^{\prime }(x),1)$, $x\in S^{1}\times D^{n-1}$.
\end{quote}

\noindent Let $[S^{1},M]^{fr}$ be the set of isotopy classes $[S^{1},F]$ of
framed circles on $M$, and write

\begin{quote}
$\gamma :$ $[S^{1},M]^{fr}\rightarrow \pi _{1}(M)$, $\gamma
([S^{1},F])=[S^{1}]$,
\end{quote}

\noindent for the projection of isotopy classes to homotopy classes. Observe
that if $\widetilde{\alpha }$ is the diffeomorphism of $S^{1}\times D^{n-1}$
in (1.1), then a framed circle $(S^{1},F)$ on $M$ determines a new one $%
(S^{1},F\widetilde{\alpha })$.

\bigskip

\noindent \textbf{Lemma 2.1.} \textsl{If }$n\geq 5$\textsl{\ the map }$%
\gamma $ \textsl{surjects. Moreover, for any framed circle }$(S^{1},F)$%
\textsl{\ on} $M$\textsl{\ we have }$\gamma ^{-1}[S^{1}]=\{[S^{1},F],[S^{1},F%
\widetilde{\alpha }]\}$\textsl{.}

\bigskip

\noindent \textbf{Proof.} With $n\geq 5$ any map $S^{1}\rightarrow M$ is
homotopic to a smooth embedding, and any two homotopic embeddings $%
S^{1}\subset M$ are isotopic (by a result of Whitney). Furthermore, since
both $S^{1}$ and $M$ are oriented, the normal disc bundle $N(S^{1})$ of an
embedding $S^{1}\subset M$ has the induced orientation, hence has a
trivialization $N(S^{1})\cong S^{1}\times D^{n-1}$, showing that $\gamma $
surjects.

As the group of trivializations of $N(S^{1})$ is $\pi _{1}(SO(n-1))=\mathbb{Z%
}_{2}$, which is generated by $\alpha $, the pre-image $\gamma ^{-1}[S^{1}]$
of any $S^{1}\subset M$ contains at most two elements: if one of them is $%
[S^{1},F]$, the other must be $[S^{1},F\widetilde{\alpha }]$.\hfill $\square 
$

\bigskip

By the proof of Lemma 2.1 we may consider a framed circle $(S^{1},F)$ as a
tubular neighborhood of the embedding $S^{1}\subset M$; project on the
factor $D^{n-1}$, and shrink the boundary to a point $\infty $, giving a
sphere $S^{n-1}$, and finally extend the map of the tubular neighborhood to $%
M$ by mapping the complement to $\infty $. We have defined a map $%
M\rightarrow S^{n-1}$ whose homotopy class is denoted by $\gamma
_{1}(S^{1},F)\in \lbrack M,S^{n-1}]$. We use $\varepsilon (S^{1},F)=0$ or $1$
to denote whether $\gamma _{1}(S^{1},F)=0$ or not.

\bigskip

\noindent \textbf{Lemma 2.2.} \textsl{Let }$M\in \mathcal{T}_{n}$\textsl{\
be }$1$\textsl{-connected with }$n\geq 5$\textsl{. Then the map }$%
\varepsilon :[S^{1},M]^{fr}\rightarrow \mathbb{Z}_{2}$\textsl{\ injects, and
satisfies the following relations:}

\textsl{i) If }$\varepsilon (S^{1},F)=1$\textsl{\ then }$\varepsilon (S^{1},F%
\widetilde{\alpha })=0$\textsl{;}

\textsl{ii) If }$w_{2}(M)=0$\textsl{\ and }$\varepsilon (S^{1},F)=0$\textsl{%
\ then }$\varepsilon (S^{1},F\widetilde{\alpha })=1$\textsl{;}

\textsl{iii) If }$w_{2}(M)\neq 0$\textsl{\ and }$\varepsilon (S^{1},F)=0$%
\textsl{\ then }$\varepsilon (S^{1},F\widetilde{\alpha })=0$\textsl{.}

\textsl{In addition,} \textsl{for any} $[S^{1},F]\in \lbrack S^{1},M]^{fr}$ 
\textsl{we can write }$M=M\#S^{n}$\textsl{\ in which the sphere }$S^{n}$ 
\textsl{has the decomposition}

\begin{enumerate}
\item[(2.1)] $S^{n}=\left\{ 
\begin{tabular}{l}
$D^{2}\times S^{n-2}\cup _{id}(S^{1}\times D^{n-1})$ \textsl{if} $%
\varepsilon (S^{1},F)=0$\textsl{;} \\ 
$D^{2}\times S^{n-2}\cup _{\tau }(S^{1}\times D^{n-1})$ \textsl{if} $%
\varepsilon (S^{1},F)=1$\textsl{,}%
\end{tabular}%
\right. $
\end{enumerate}

\noindent \textsl{such that the second component }$S^{1}\times D^{n-1}$%
\textsl{\ is identical to }$F$\textsl{.}

\bigskip

\noindent \textbf{Proof. }The first assertion, together with the properties
i)-iii), have been shown by Goldstein and Linger in \cite[Theorem 1]{GL}.

Let $S^{1}\subset $\textsl{\ }$M$ be an oriented circle on $M$. Since $M$%
\textsl{\ }is $1$-connected with dimension\textsl{\ }$n\geq 5$ there is an
embedded disc $D^{2}\subset M$ with $\partial D^{2}=S^{1}\subset M$. Taking
a tubular neighborhood $D^{n}$ of the disc $D^{2}$ one gets a frame $F$ of
the circle $S^{1}$ such that $\func{Im}F\subset \overset{\circ }{D^{n}}$.
That is, the manifold $M$ has a decomposition $M=M\#S^{n}$ so that $%
(S^{1},F) $ is a framed circle on the sphere $S^{n}$.

For $M=S^{n}$ we have $[M,S^{n-1}]=\mathbb{Z}_{2}$, while the map $%
\varepsilon $ is one to one, implying that $S^{n}$ possesses the two
decompositions stated in (2.1).\hfill $\square $

\bigskip

\noindent \textbf{Corollary 2.3. }\textsl{If }$M\in \mathcal{T}_{n}$\textsl{%
\ is }$1$\textsl{-connected with }$n\geq 5$\textsl{, then }$[S^{1},M]^{fr}=%
\mathbb{Z}_{2}$\textsl{\ or }$0$\textsl{, where }$[S^{1},M]^{fr}=\mathbb{Z}%
_{2}$\textsl{\ happens if and only if }$w_{2}(M)=0$\textsl{.}\hfill $\square 
$

\bigskip

For a $N\in \mathcal{T}_{n-1}$ we furnish the product $S^{1}\times N$ with
the induced orientation, and let $(D^{n-1},0)\subset (N,x_{0})$ be an
orientation preserving embedded disc centered at $x_{0}$. Then the natural
inclusion $F_{0}:S^{1}\times D^{n-1}\subset S^{1}\times N$ gives rise to the 
\textsl{canonical frame }of the oriented circle $S^{1}\times \{x_{0}\}$ on $%
S^{1}\times N$, for which $\gamma (S^{1}\times \{x_{0}\},F_{0})$ is the
generator of the first summand of $\pi _{1}(S^{1}\times N)=\mathbb{Z}\oplus
\pi _{1}(N)$ that corresponds to the orientation on $S^{1}\times \{x_{0}\}$,
while $\gamma _{1}(S^{1}\times \{x_{0}\},F_{0})$ is the homotopy class of
the composition

\begin{quote}
$c:S^{1}\times N\rightarrow N\rightarrow N/(N-\overset{\circ }{D_{n}}%
)=S^{n-1}$.
\end{quote}

\noindent For a frame $F$ of the circle $S^{1}\times \{x_{0}\}$ on $%
S^{1}\times N$ we write $\delta (S^{1}\times \{x_{0}\},F)=0$ or $1$ to
denote the statement that $\gamma _{1}(S^{1}\times \{x_{0}\},F)=[c]$ or not.
From Lemma 2.1 we get

\bigskip

\noindent \textbf{Lemma 2.4.} \textsl{If }$n\geq 4$\textsl{\ and} \textsl{if 
}$F$\textsl{\ is a frame of the oriented circle }$S^{1}\times \{x_{0}\}$%
\textsl{\ on} $S^{1}\times N$\textsl{,} \textsl{then }$\delta (S^{1}\times
\{x_{0}\},F)=1$\textsl{\ implies} \textsl{that} $\delta (S^{1}\times
\{x_{0}\},F\widetilde{\alpha })=0.$\hfill $\square $

\bigskip

Now, let $(S^{1},F)$ be a framed circle on a $1$-connected manifold $M\in 
\mathcal{T}_{n}$ with $n\geq 5$, and let $F^{\prime }$ be one of the two
possible frames of the oriented circle $S^{1}\times \{x_{0}\}$\ on $%
S^{1}\times N$, $N\in \mathcal{T}_{n-1}$. \textsl{The tunnel sum} $\circ $ of%
\textsl{\ }$M$ and $S^{1}\times N$\textsl{\ }along the framed circles $%
(S^{1},F)$ and $(S^{1}\times \{x_{0}\},F^{\prime })$ is the adjoint manifold

\begin{enumerate}
\item[(2.2)] $(M,F)\circ (S^{1}\times N,F^{\prime }):=(M-\overset{\circ }{%
\func{Im}F})\cup _{f\circ f^{^{\prime }-1}}(S^{1}\times N-\overset{\circ }{%
\func{Im}F^{\prime }})$.
\end{enumerate}

\noindent where

\begin{quote}
$f:=F\mid S^{1}\times S^{n-2}$, $f^{\prime }:=F^{\prime }\mid S^{1}\times
S^{n-2}$.
\end{quote}

\noindent Useful properties of the operation $\circ $ are collected in the
following lemma.

\bigskip

\noindent \textbf{Lemma 2.5.} \textsl{The tunnel sum operation }$\circ $ 
\textsl{satisfies the following relations}

\begin{enumerate}
\item[(2.3)] $(M,F)\circ (S^{1}\times N,F^{\prime })=(M,F\widetilde{\alpha }%
)\circ (S^{1}\times N,F^{\prime }\widetilde{\alpha })$;

\item[(2.4)] $(M,F)\circ (S^{1}\times N,F^{\prime })=\left\{ 
\begin{tabular}{l}
$M\#\Sigma _{0}N$ \textsl{if} $\varepsilon (S^{1},F)=\delta (S^{1},F^{\prime
})$\textsl{;} \\ 
$M\#\Sigma _{1}N$ \textsl{if} $\varepsilon (S^{1},F)\neq \delta
(S^{1},F^{\prime })$\textsl{.}%
\end{tabular}%
\right. $
\end{enumerate}

\noindent \textsl{In addition,}

\begin{enumerate}
\item[(2.5)] $M\#\Sigma _{0}N=M\#\Sigma _{1}N$\textsl{\ if }$w_{2}(M)\neq 0$%
\textsl{.}
\end{enumerate}

\noindent \textbf{Proof.} Since the generator $\alpha \in \pi _{1}(SO(n-1))$
is of order $2$, the relation (2.3) follows from the definition of the
tunnel sum in (2.2).

Since $M$ is $1$-connected, we can assume by Lemma 2.2 that $M=S^{n}$, and
take a decomposition of $S^{n}$ relative to the index $\varepsilon (S^{1},F)$
as in (2.1). If $\varepsilon (S^{1},F)=\delta (S^{1},F^{\prime })=0$ (resp.
if $\varepsilon (S^{1},F)=1$ but $\delta (S^{1},F^{\prime })=0$) we get
(2.4) directly from the definition of $\Sigma _{i}:$

\begin{quote}
$(S^{n},F)\circ (S^{1}\times N,F^{\prime })=D^{{\small 2}}\times S^{{\small %
n-1}}\cup _{id}(S^{{\small 1}}\times (N-\overset{\circ }{D^{n-1}}))=\Sigma
_{0}N$

(resp. $=D^{{\small 2}}\times S^{{\small n-1}}\cup _{\tau }(S^{{\small 1}%
}\times (N-\overset{\circ }{D^{n-1}})=\Sigma _{1}N$).
\end{quote}

\noindent For the remaining case $\varepsilon (S^{1},F)=\delta
(S^{1},F^{\prime })=1$ (resp. $\varepsilon (S^{1},F)=0$ but $\delta
(S^{1},F^{\prime })=1$), we can compose both $F$ and $F^{\prime }$ with the
diffeomorphism $\widetilde{\alpha }$ of $S^{{\small 1}}\times D^{n-1}$ by
(2.3), then apply i) of Lemma 2.2 and Lemma 2.4 to reduce the current case
to the previous one.

Finally, the relation (2.5) has been shown by iii) of Lemma 2.2.\hfill $%
\square $

\bigskip

We are ready to show Theorem B stated in Section 1.

\bigskip

\noindent \textbf{Proof of Theorem B. }Let $p:E\rightarrow B$ an oriented
circle bundle over a manifold $B\in \mathcal{T}_{n-1}$, where\textbf{\ }the
total space\textbf{\ }$E$\ is $1$-connected with $\dim E=n\geq 5$. Take an
orientation preserving embedding $i:D^{n-1}\subset B$ and set $i(0)=b\in B$.
Then the inclusion $F_{b}:E\mid D^{n-1}\subset E$ defines a frame $%
(p^{-1}(b),F_{b})$ of the oriented circle $p^{-1}(b)\subset E$, to be called
a \textsl{canonical frame of }$p^{-1}(b)$. Recall from \cite[Theorem 8]{GL}
that the function $\varepsilon :B\rightarrow \{0,1\}$ by $\varepsilon
(b):=\varepsilon (p^{-1}(b),F_{b})$ is constant, and satisfies that

\begin{enumerate}
\item[(2.6)] $\varepsilon (p^{-1}(b),F_{b})=0$ or $1$ in accordance to $%
w_{2}(B)\neq 0$ or $w_{2}(B)=0$.
\end{enumerate}

On the other hand, let $p^{\prime }:S^{1}\times N\rightarrow N$\ be the
trivial $S^{1}$-bundle on $N\in \mathcal{T}_{n-1}$. Take an orientation
preserving embedding $i:D^{n-1}\subset N$ and set $x_{0}:=i(0)\in N$. Since
the inclusion $F_{x_{0}}:p^{\prime }\mid D^{n-1}\subset S^{1}\times N$ is
precisely the canonical frame $(p^{\prime -1}(x_{0}),F_{x_{0}})$ of the
oriented circle $p^{\prime -1}(x_{0})=S^{1}\times x_{0}$ on $S^{1}\times N$
we get by the definition of the index $\delta $ that

\begin{enumerate}
\item[(2.7)] $\delta (p^{\prime -1}(x_{0}),F_{x_{0}})=0$.
\end{enumerate}

Consider now the projection $f_{N}:B\#N\rightarrow B$ ( $f_{N}(N-\overset{%
\circ }{D}^{n-1})=b$), and let $p_{N}:E_{N}\rightarrow $ $B\#N$ be the
induced bundle $f_{N}^{\ast }E$. Then the total space $E_{N}$ admits the
decomposition into the tunnel sum

\begin{quote}
$E_{N}=(E,F_{b})\circ (S^{1}\times N,F_{x_{0}})$,
\end{quote}

\noindent where $b\in B$, $x_{0}\in N$, and where both $(p^{-1}(b),F_{b})$
and $(S^{1}\times \{x_{0}\},F_{x_{0}})$ are canonical. In view of (2.6) and
(2.7), the desired formulae (1.4)\ and (1.5) of Theorem B have been shown by
the relations (2.4) and (2.5) in Lemma 2.5, respectively.\hfill $\square $

\section{Relationship between $\Sigma _{1}N$\ and $\Sigma _{0}N$}

This section proves Theorem A, and illustrates related applications and
examples.

\bigskip

\noindent \textbf{Proof of Theorem A.} Let $N\in \mathcal{T}_{n}$\ be a
manifold with base point $x_{0}\in N$, and let $F_{0}$ be the canonical frame%
\textsl{\ }of the oriented circle $S^{1}\times \{x_{0}\}$ on $S^{1}\times N$%
. A map $\beta :S^{1}\rightarrow Diff(N,x_{0})$ gives rise to a
diffeomorphism of $S^{1}\times N$

\begin{quote}
$\widetilde{\beta }:S^{1}\times N\rightarrow S^{1}\times N$ by $\widetilde{%
\beta }(z,x)=(z,\beta (z)\cdot x)$
\end{quote}

\noindent that fixes the circle $S^{1}\times x_{0}$ pointwisely. In
particular, in addition to $F_{0}$ the composition $\widetilde{\beta }\circ
F_{0}$ is a second frame of the circle $S^{1}\times \{x_{0}\}$.

Assume that the manifold $N$ is $\Sigma $-stable. Then there exists a $\beta 
$ such that $e_{N}\circ \beta =\alpha $, implying $\widetilde{\beta }\circ
F_{0}=F_{0}\widetilde{\alpha }$. Thus, surgeries along $S^{1}\times x_{0}$
with respect to the two frames $F_{0}$ and $\widetilde{\beta }\circ F_{0}$,
we obtain from $\widetilde{\beta }$ a desired diffeomorphism $\Sigma
_{0}N\rightarrow \Sigma _{1}N$.\hfill $\square $

$\bigskip $

For two manifolds $N\in \mathcal{T}_{n}$ and $M\in \mathcal{T}_{m}$ with
base points $x_{0}\in N$ and $y_{0}\in M$ the linear representation $%
e_{N\times M}$ fits into the commutative diagram

\begin{quote}
$%
\begin{array}{ccc}
Diff(N,\{x_{0}\}) & \overset{e_{N}}{\rightarrow } & GL(n) \\ 
i\downarrow &  & j\downarrow \\ 
Diff(N\times M,\{x_{0},y_{0}\}) & \overset{e_{N\times M}}{\rightarrow } & 
GL(n+m)%
\end{array}%
$,
\end{quote}

\noindent where $i(h)=h\times id_{M}$ and $j(g)=g\oplus I_{m}$. By the fact
that $j$ is an $(n-1)$-homotopy equivalence we get

\bigskip

\noindent \textbf{Corollary 3.1.}\textsl{\ If }$N\in \mathcal{T}_{n}$\textsl{%
\ with }$n\geq 3$ \textsl{is }$\Sigma $\textsl{-stable, then the product }$%
N\times M$ \textsl{is also} $\Sigma $\textsl{-stable. }

\textsl{In particular, for any }$M\in \mathcal{T}_{m}$\textsl{,} $\Sigma
_{0}(N\times M)=\Sigma _{1}(N\times M)$\textsl{.}\hfill $\square $

\bigskip

For a manifold $N\in T_{n}$\ with an orientation preserving embedding $%
(D^{n},0)\subset (N,x_{0})$ set $W_{N}:=D^{2}\times (N-\overset{\circ }{D^{n}%
})$. In addition to (1.2) we have (e.g. \cite{Su})

\begin{enumerate}
\item[(3.1)] $\Sigma _{0}N=\partial W_{N}$.
\end{enumerate}

\noindent This relation, together with Theorem A, is useful to calculate the
diffeomorphism types of some $\Sigma _{i}N$. We present such examples.

\bigskip

\noindent \textbf{Proposition 3.2.} \textsl{The} \textsl{sphere }$S^{n}$%
\textsl{\ is} $\Sigma $\textsl{-stable with }$\Sigma _{0}S^{n}=\Sigma
_{1}S^{n}=S^{n+1}$.

\textsl{In particular,} \textsl{for any }$p\leq q$\textsl{\ with }$q\geq 3$ 
\textsl{we have}

\begin{enumerate}
\item[(3.2)] $\Sigma _{0}(S^{p}\times S^{q})=\Sigma _{1}(S^{p}\times
S^{q})=S^{p}\times S^{q+1}\#S^{p+1}\times S^{q}$.
\end{enumerate}

\noindent \textbf{Proof. }Let\textbf{\ }$\{e_{1},\cdots ,e_{n+1}\}$ be the
standard orthonormal basis of $\mathbb{R}^{n+1}$ and take $e_{n+1}\in S^{n}$
as the base point. The\textbf{\ }generator\textbf{\ }$\alpha \in \pi
_{1}(SO(n))$ gives rise to the loop on $Diff(S^{n},x_{0})$

\begin{quote}
$\beta :S^{1}\rightarrow Diff(S^{n},x_{0})$, $\beta (z)(x)=\alpha (z)x$, $%
z\in S^{1}$, $x\in S^{n}$,
\end{quote}

\noindent that satisfies $e_{S^{n}}(\beta )=\alpha $, implying that $S^{n}$%
\textsl{\ }is $\Sigma $-stable. We get $\Sigma _{0}S^{n}=\Sigma _{1}S^{n}$
by Theorem A. Moreover, from $W_{S^{n}}=D^{n+2}$ we get by (3.1) that

\begin{quote}
$\Sigma _{0}S^{n}=\partial W_{S^{n}}=S^{n+1}$,
\end{quote}

\noindent showing the first statement of the proposition.

For (3.2) we get from (3.1) that (e.g. \cite[Lemma 1.3]{Su}, \cite[Lemma 2]%
{G})

\begin{quote}
$\Sigma _{0}(S^{p}\times S^{q})=\partial W_{S^{p}\times S^{q}}=S^{p}\times
S^{q+1}\#S^{p+1}\times S^{q}$.
\end{quote}

\noindent Since $S^{p}\times S^{q}$ is $\Sigma $\textsl{-}stable by
Corollary 3.1, we get (3.2) by Theorem A.\hfill $\square $

\bigskip

Let $\mathbb{F}$ denote either the field $\mathbb{R}$ of reals, the field $%
\mathbb{C}$ of complexes, or the algebra $\mathbb{H}$ of quaternions. Let $%
\mathbb{F}P^{n}$ be the $n$-dimensional projective space of the $1$%
-dimensional $\mathbb{F}$-subspaces of the vector space $\mathbb{F}^{n+1}$.
The canonical Hopf line bundle over $\mathbb{F}P^{n}$ is

\begin{quote}
$\lambda _{\mathbb{F}}^{n}:=\{(x,v)\in \mathbb{F}P^{n}\times \mathbb{F}%
^{n+1}\mid v\in x\}$.
\end{quote}

\noindent In accordance to $\mathbb{F}=\mathbb{R},\mathbb{C}$ or $\mathbb{H}$
the real reduction of $\lambda _{\mathbb{F}}^{n}$ is a real Euclidean vector
bundle over $\mathbb{F}P^{n}$ with dimension $1$, $2$ or $4$, respectively.
For an Euclidean vector bundle $\xi $ over a CW-complex $B$ write $D(\xi )$
and $S(\xi )$ to denote the total spaces of the associated disc-bundle and
sphere-bundle, respectively.

\bigskip

\noindent \textbf{Proposition 3.3. }$\Sigma _{0}\mathbb{F}P^{n}=S(\lambda _{%
\mathbb{F}}^{n-1}\oplus \varepsilon ^{2})$\textsl{, where }$\varepsilon ^{2}$%
\textsl{\ denotes the }$2$\textsl{-dimensional trivial bundle over }$\mathbb{%
F}P^{n-1}$\textsl{. }

\textsl{In particular, the homotopy groups of the suspension }$\Sigma _{0}%
\mathbb{F}P^{n}$ \textsl{are}

\begin{enumerate}
\item[(3.3)] $\pi _{r}(\Sigma _{0}\mathbb{F}P^{n})=$\textsl{\ }$\pi _{r}(%
\mathbb{F}P^{n-1})\oplus \pi _{r}(S^{d})$\textsl{,}
\end{enumerate}

\noindent \textsl{where }$d=2$\textsl{, }$3$\textsl{\ or }$5$\textsl{\ in
accordance to }$\mathbb{F}=\mathbb{R},\mathbb{C}$\textsl{\ or }$\mathbb{H}$%
\textsl{.}

\bigskip

\noindent \textbf{Proof.} Set $r:=\dim _{\mathbb{R}}\mathbb{F}$. In view of $%
\partial D(\lambda _{\mathbb{F}}^{n-1})=S^{rn-1}$ the manifold $\mathbb{F}%
P^{n}$ has the decomposition $\mathbb{F}P^{n}=D(\lambda _{\mathbb{F}%
}^{n-1})\cup _{id}D^{rn}$, implying

\begin{quote}
$W_{\mathbb{F}P^{n}}=D^{2}\times (\mathbb{F}P^{n}-\overset{\circ }{D^{rn}}%
)=D^{2}\times D(\lambda _{\mathbb{F}}^{n-1})=D(\lambda _{\mathbb{F}%
}^{n-1}\oplus \varepsilon ^{2})$.
\end{quote}

\noindent One obtains by (3.1) that

\begin{quote}
$\Sigma _{0}\mathbb{F}P^{n}=\partial W_{\mathbb{F}P^{n}}=S(\lambda _{\mathbb{%
F}}^{n-1}\oplus \varepsilon ^{2})$.
\end{quote}

\noindent Since $\Sigma _{0}\mathbb{F}P^{n}$ is the total space of a $d$%
-dimensional sphere bundle over $\mathbb{F}P^{n-1}$ which has a section, we
obtain (3.3) by the homotopy exact sequence of the fiber bundle $\Sigma _{0}%
\mathbb{F}P^{n}\rightarrow \mathbb{F}P^{n-1}$.\hfill $\square $

\bigskip

\noindent \textbf{Example 3.4.} Let $\Omega _{n}^{Spin}$ be the $n$%
-dimensional Spin bordism group. Then $\Sigma _{0}\mathbb{H}P^{2}=0$ in $%
\Omega _{9}^{Spin}$ by $\Sigma _{0}\mathbb{H}P^{2}=\partial W_{\mathbb{H}%
P^{2}}$. On the other hand, there exists a surgery along the circle $%
S^{1}\times \{x_{0}\}\subset S^{1}\times \mathbb{H}P^{2}$ such that the
resulting manifold $M$ satisfies $M\neq 0$ in $\Omega _{9}^{Spin}$ (e.g. 
\cite[p.258]{AB}, \cite{DY}). That is

\begin{quote}
$M\neq \Sigma _{0}\mathbb{H}P^{2}$ but $M=\Sigma _{1}\mathbb{H}P^{2}$.
\end{quote}

\noindent In particular, the manifold $\mathbb{H}P^{2}$ fails to be $\Sigma $%
-stable by Theorem A.\hfill $\square $

\bigskip

\noindent \textbf{Example 3.5. }Formula (1.10) in Section 1 suggests that
Theorem B may be applicable to derive the homotopy groups of certain
manifolds. Combining (1.6) with (3.3) we get further evidences:

\begin{quote}
$\pi _{r}(\mathbb{C}P^{2n}\#\mathbb{R}P^{4n})=\left\{ 
\begin{array}{c}
\mathbb{Z}_{2}\text{ if }r=1\text{,} \\ 
\pi _{r}(S^{4n-1})\oplus \pi _{r}(S^{2})\text{ if }r\neq 1%
\end{array}%
\right. $;

$\pi _{r}(\mathbb{C}P^{2n}\#\mathbb{C}P^{2n})=\left\{ 
\begin{array}{c}
\mathbb{Z}\oplus \mathbb{Z}\text{ if }r=2 \\ 
\pi _{r}(S^{4n-1})\oplus \pi _{r}(S^{3})\text{ if }r\neq 2%
\end{array}%
\right. $;

$\pi _{r}(\mathbb{C}P^{2n}\#\mathbb{H}P^{n})=\left\{ 
\begin{array}{c}
\mathbb{Z}\text{ if }r=2,4 \\ 
\pi _{r}(S^{4n-1})\oplus \pi _{r}(S^{5})\text{ if }r\neq 2,4%
\end{array}%
\right. $.\hfill $\square $
\end{quote}

\bigskip

\noindent \textbf{Example 3.6.} For $n=2$ we have $\pi _{1}(SO(2))=\mathbb{Z}
$ with generator $\alpha $. Therefore, the operation $\Sigma _{i}$ applies
also to the orientable surfaces to yield the correspondences $\Sigma _{i}:%
\mathcal{T}_{2}\rightarrow \mathcal{T}_{3}$, $i=0,1$.

Let $M_{g}\in \mathcal{T}_{2}$ be the orientable surface with genus $g\geq 1$%
, and let $D^{2}\subset M_{g}$ be an embedded disc. Then the fundamental
group $\pi _{1}(M_{g}-\overset{\circ }{D^{2}})$ is the free group generated
by the standard elements

\begin{quote}
$a_{1},b_{1},\cdots a_{g},b_{g}\in \pi _{1}(M_{g}-\overset{\circ }{D^{2}})$.
\end{quote}

\noindent Applying the Van Kampen Theorem to (1.2) one gets the
presentations of the fundamental groups $\pi _{1}(\Sigma _{i}M_{g})$ by the
generators and relations \cite{F}

\begin{quote}
$\pi _{1}(\Sigma _{0}M_{g})=\left\langle a_{1},b_{1},\cdots
a_{g},b_{g}\right\rangle $;

$\pi _{1}(\Sigma _{1}M_{g})=\left\langle a_{1},b_{1},\cdots
a_{g},b_{g},z\mid a_{i}z=za_{i},b_{i}z=zb_{i},\right. $

$\qquad \qquad \qquad \left. z^{-1}=a_{1}b_{1}a_{1}^{-1}b_{1}^{-1}\cdots
a_{n}b_{n}a_{n}^{-1}b_{n}^{-1}\right\rangle $.
\end{quote}

\noindent Precisely, $\Sigma _{0}M_{g}=S^{1}\times S^{2}\#\cdots
\#S^{1}\times S^{2}$ ($2g$-copies), while $\Sigma _{1}M_{g}$ is the total
space of the oriented circle bundle over $M_{g}$ whose Euler class is the
orientation class on $M_{g}$ (e.g. Freedman \cite{F}). This shows that $%
\Sigma _{0}M_{g}\neq \Sigma _{1}M_{g}$.\hfill $\square $

\section{Properties of the suspensions $\Sigma _{i}$}

This section develops general properties of the suspensions $\Sigma _{i}$,
required to show Theorems C and D.

\bigskip

\noindent \textbf{Proposition 4.1.}\textsl{\ For a set of }$1$\textsl{%
-connected manifolds }$N_{1},\cdots ,N_{k}\in \mathcal{T}_{n-1}$\textsl{\
with }$\QTR{sl}{n\geq 5}$ \textsl{we have}

\begin{enumerate}
\item[(4.1)] $\Sigma _{i}(N_{1}\#\cdots \#N_{k})=\Sigma _{i}N_{1}\#\cdots
\#\Sigma _{i}N_{k}$, $i=0,1$.
\end{enumerate}

\noindent \textbf{Proof. }It suffices to show (4.1) for the case of $i=1$.
The remaining case $i=0$ follows from the same idea.

By the definition of the space $E_{N}$ in (1.3) we have $E_{N%
\#M}=(E_{N})_{M} $. In addition, the Van-Kampen Theorem \cite[p.43]{Ha}
implies that, if $N\in \mathcal{T}_{n-1}$ is $1$-connected, then

\begin{quote}
$\pi _{1}(\Sigma _{1}N)=\pi _{1}(N)=0$, $i=0,1$.
\end{quote}

\noindent Since $\Sigma _{1}N_{1},\cdots ,\Sigma _{1}N_{k}\in \mathcal{T}%
_{n-1}$ are all $1$-connected, Theorem B can be applied repeatedly to show
that, if $w_{2}(B)=0$, then

\begin{quote}
$E_{N_{1}\#\cdots \#N_{k}}=E\#\Sigma _{1}N_{1}\#\cdots \#\Sigma _{1}N_{k}$.
\end{quote}

\noindent In particular, taking $p$ to be the circle bundle $S^{3}\times
S^{n-2}\rightarrow S^{2}\times S^{n-2}$ in Example 1.5, we get

\begin{quote}
$S^{3}\times S^{n-2}\#\Sigma _{1}(N_{1}\#\cdots \#N_{k})=S^{3}\times
S^{n-2}\#\Sigma _{1}N_{1}\#\cdots \#\Sigma _{1}N_{k}$.
\end{quote}

\noindent Since the surgery on the $3$-sphere $S^{3}\times \{\ast \}\subset $
$S^{3}\times S^{n-3}$ with trivial framing yields $S^{n}$, applying surgery
on both sides to kill $S^{3}\times \{\ast \}$ shows (4.1).\hfill $\square $

\bigskip

A manifold $N$ is a \textsl{homology }$n$\textsl{-sphere} if its integral
homology is isomorphic to that of the $n$-sphere $S^{n}$. The suspension $%
\Sigma _{i}$ enables one to produce higher dimensional homology spheres from
the lower dimensional ones. For a homology $n$-sphere $N$ consider the
circle bundle given by Example 1.5:

\begin{quote}
$p_{N}:S^{3}\times S^{n-2}\#\Sigma _{i}N\rightarrow S^{2}\widetilde{\times }%
S^{n-2}\#N$ (resp. $\rightarrow S^{2}\times S^{n-2}\#N$).
\end{quote}

\noindent By the Gysin sequence of $p_{N}$ \cite[p.143]{Mil2} with integer
coefficients we get

\begin{quote}
$H^{\ast }(S^{3}\times S^{n-2}\#\Sigma _{i}N)=H^{\ast }(S^{3}\times S^{n-2})$%
.
\end{quote}

\noindent This shows that

\bigskip

\noindent \textbf{Proposition 4.2.} \textsl{If }$N$\textsl{\ is a homology }$%
n$\textsl{-sphere with }$n\geq 4$\textsl{, then both }$\Sigma _{i}N$\textsl{%
, }$i=0,1$\textsl{,\ are homology }$(n+1)$\textsl{-sphere.}\hfill $\square $

\bigskip

We conclude this section by\ expressing the topological invariants (i.e. the
integral homology, cohomology, and characteristic classes) of the
suspensions $\Sigma _{i}N$\ in term of that of $N$. Recall that \textsl{the
reduced suspension} $SX$ of a CW-complex $X$ is the quotient space $%
S^{1}\times X/(S^{1}\vee X)$ \cite[p.223]{Ha}. Let $\widetilde{H_{\ast }}$
(resp. $\widetilde{H^{\ast }}$) denotes the reduced homology (resp.
cohomology) of topological spaces.

\bigskip

\noindent \textbf{Proposition 4.3.} \textsl{The homology and cohomology of
the suspensions }$\Sigma _{i}N$\textsl{\ have the following presentations}

\begin{enumerate}
\item[(4.2)] $H_{\ast }(\Sigma _{i}N)=\widetilde{H_{\ast }}(SN)\oplus
H_{\ast }(N-\left\{ x_{0}\right\} )$;

\item[(4.3)] $H^{\ast }(\Sigma _{i}N)=\widetilde{H^{\ast }}(SN)\oplus
H^{\ast }(N-\left\{ x_{0}\right\} )$.
\end{enumerate}

\textsl{In addition, with respect to the decomposition of }$H^{\ast }(\Sigma
_{i}N)$ \textsl{in (4.3),} \textsl{the Stiefel-Whitney classes }$w_{r}$%
\textsl{\ and Pontryagin classes }$p_{r}$\textsl{\ of }$\Sigma _{i}N$\textsl{%
\ are given respectively by the formulae}

\begin{enumerate}
\item[(4.4)] $w_{r}(\Sigma _{i}N)\equiv 0\oplus w_{r}(N-\left\{
x_{0}\right\} )$ $\func{mod}$ $\mathbb{Z}_{2}$\textsl{; }

$p_{r}(\Sigma _{i}N)=0\oplus p_{r}(N-\left\{ x_{0}\right\} )$\textsl{.}
\end{enumerate}

\noindent \textbf{Proof. }It suffices to show formulae (4.3) and (4.4).%
\textbf{\ }In view of (1.2) we set $A:=S^{1}\times (N-\overset{\circ }{D^{n}}%
)$ and $B:=D^{2}\times S^{n-1}$ to get

\begin{quote}
$\Sigma _{i}N=A\cup B$, $A\cap B=S^{1}\times S^{n-1}$,
\end{quote}

\noindent as well as the Mayer-Vietoris sequence

\begin{center}
$\cdots \rightarrow H^{r-1}(A\cap B)\overset{\delta }{\rightarrow }%
H^{r}(\Sigma _{i}N)\overset{i_{1}^{\ast }\oplus i_{2}^{\ast }}{\rightarrow }%
H^{r}(A)\oplus H^{r}(B)\overset{j_{1}^{\ast }-j_{2}^{\ast }}{\rightarrow }%
H^{r}(A\cap B)\rightarrow \cdots $.
\end{center}

\noindent For a $M\in \mathcal{T}_{n}$ let $\omega _{M}\in H^{n}(M)$ be the
orientation class. Then the graded group $H^{\ast }(A\cap B)$ is torsion
free with basis 

\begin{quote}
$\{1,\omega _{S^{1}}\times 1,1\times \omega _{S^{n-1}},\omega _{S^{1}\times
S^{n-1}}\}$ 
\end{quote}

\noindent where $\omega _{S^{1}}\times 1\in \func{Im}j_{1}^{\ast }$, $%
1\times \omega _{S^{n-1}}\in \func{Im}j_{2}^{\ast }$. Moreover, the
Mayer-Vietoris sequence can be summarized as the short exact sequence

\begin{enumerate}
\item[(4.5)] $0\rightarrow \mathbb{Z}\{\omega _{S^{1}\times S^{n-1}}\}%
\overset{\delta }{\rightarrow }H^{\ast }(\Sigma _{i}N)\overset{i_{1}^{\ast }}%
{\rightarrow }\frac{H^{\ast }(S^{1}\times (N-\left\{ x_{0}\right\} ))}{%
\mathbb{Z}\{\omega _{S^{1}}\times 1\}}\rightarrow 0$,
\end{enumerate}

\noindent where

a) $\mathbb{Z}\{\omega \}$ denotes the cyclic group with generator $\omega $;

b) $N-\overset{\circ }{D^{n}}$ has been replaced by $N-\left\{ x_{0}\right\} 
$, and

c) the map $\delta $ is an isomorphism onto $H^{n+1}(\Sigma _{i}N)=\mathbb{Z}
$.

\noindent In particular, the sequence (4.5) is splittable to yields (4.3),
under the convention that the orientation class $\omega _{\Sigma
_{i}N}=\delta (\omega _{S^{1}\times S^{n-1}})$ agrees with the top degree
generator $s\wedge \omega _{N}$ of $\widetilde{H^{n+1}}(SN)=\mathbb{Z}$.

Finally, let $TM$ be the tangent bundle of a smooth manifold $M$. Then the
inclusion $i_{1}:$ $S^{1}\times (N-\overset{\circ }{D^{n}})\rightarrow
\Sigma _{i}N$ satisfies that

\begin{quote}
$i_{1}^{\ast }T(\Sigma _{i}N)=\varepsilon ^{1}\oplus T(N-\left\{
x_{0}\right\} )$,
\end{quote}

\noindent where $\varepsilon ^{1}$ is the $1$-dimensional trivial bundle on $%
S^{1}\times (N-\overset{\circ }{D^{n}})$ consisting of the tangent
directions of the first factor $S^{1}$. The formulae (4.4) follows now from
the naturality of the characteristic classes, together with the
decomposition (4.3).\hfill $\square $

\section{Circle actions on the 1-connected 6-manifolds}

In \cite{GL1,L,GL} R. Goldstein and L. Lininger initiated the project to
classify those 1-connected 6-manifolds that admit regular circle actions. In
this section we complete the project by showing Theorem C, which makes the
main results of \cite{GL1,GL} precise.

Let $A$ be an abelian group. An element $\omega \in A$ is called \textsl{%
primitive} if the cyclic subgroup generated by $\omega $ is a direct summand
of $A$. By the homotopy exact sequence of fiber bundles one gets (e.g. \cite%
{DL})

\bigskip

\noindent \textbf{Lemma 5.1. }\textsl{Let }$M$\textsl{\ be an }$1$\textsl{%
-connected }$6$\textsl{-manifolds that has a regular circle action with
quotient map }$\pi :$\textsl{\ }$M\rightarrow M/S^{1}$\textsl{. Then}

\textsl{i) }$M/S^{1}$\textsl{\ is an }$1$\textsl{-connected }$5$\textsl{%
-manifolds;}

\textsl{ii) The Euler class }$\omega \in H^{2}(M/S^{1})$\textsl{\ of }$\pi $%
\textsl{\ is primitive.}\hfill $\square $

\bigskip

Let $\gamma $ be the Hopf complex line bundle over $S^{2}$, and consider the
two handlebodies of dimension $5$

\begin{quote}
$A:=S^{2}\times D^{3}$, $B:=D(\gamma \oplus \varepsilon )$,
\end{quote}

\noindent where $\varepsilon $ is the $1$-dimensional trivial bundle over $%
S^{2}$, and $D(\eta )$ denotes the disc bundle of an Euclidean vector bundle 
$\eta $. For a connected manifold $V$ with non-empty boundary $\partial V$
write $V\cdot V$ to denote \textsl{the boundary connected sum} of two copies
of $V$. As examples, one gets from $\partial A=S^{2}\times S^{2}$ and $%
\partial B=\mathbb{C}P^{2}\#\overline{\mathbb{C}P}^{2}$ that

\begin{quote}
$\partial (A\cdot A)=S^{2}\times S^{2}\#S^{2}\times S^{2}$, $\partial
(B\cdot B)=\mathbb{C}P^{2}\#\overline{\mathbb{C}P}^{2}\#\mathbb{C}P^{2}\#%
\overline{\mathbb{C}P}^{2}$.
\end{quote}

Denote by $\{a,b\}$, $\{a_{1},b_{1},a_{2},b_{2}\}$ and $%
\{x_{1},y_{1},x_{2},y_{2}\}$ the canonical bases of the second cohomologies
of the boundaries $\partial B$, $\partial (A\cdot A)$ and $\partial (B\cdot
B)$, respectively. With these convention we introduced in the table below a
collection of five $1$-connected $5$-manifolds

\begin{enumerate}
\item[(5.1)] 
\begin{tabular}{l|l|l|l|l}
\hline
$S^{2}\times S^{3}$ & $S^{2}\widetilde{\times }S^{3}$ & $W$ & $M_{k}$ & $%
X_{2^{i}}$ \\ \hline
$A\cup _{id}A$ & $B\cup _{id}B$ & $B\cup _{\tau }B$ & $A\cdot A\cup
_{f_{k}}A\cdot A$ & $B\cdot B\cup _{g_{k}}B\cdot B$ \\ \hline
\end{tabular}
\end{enumerate}

\noindent where $\tau \in Diff(\partial B)$, $f_{k}\in Diff(\partial (A\cdot
A))$, and $g_{k}\in Diff(\partial (B\cdot B))$ are the diffeomorphisms whose
actions on the second homologies are (e.g. \cite{Wall3})

\begin{quote}
i) $(\tau _{\ast }(a),\tau _{\ast }(b))=(a,-b)$;

ii) $(f_{k\ast }(a_{1}),f_{k\ast }(b_{1}),f_{k\ast }(a_{2}),f_{k\ast
}(b_{2}))=(a_{1},b_{1},a_{2},b_{2})\cdot C_{k}$,

iii) $(g_{k\ast }(x_{1}),g_{k\ast }(y_{1}),g_{k\ast }(x_{2}),g_{k\ast
}(y_{2}))=(x_{1},y_{1},x_{2},y_{2})\cdot L_{k}$,
\end{quote}

\noindent and where $C_{k}$ and $L_{k}$ denote, respectively, the $4\times 4$
matrices

\begin{quote}
$C_{k}=\left( 
\begin{array}{cccc}
1 & 0 & 0 & 0 \\ 
0 & 1 & k & 0 \\ 
0 & 0 & 1 & 0 \\ 
-k & 0 & 0 & 1%
\end{array}%
\right) $ and $L_{k}=\left( 
\begin{array}{cccc}
1 & 2^{k} & 2^{k} & 0 \\ 
2^{k} & 1 & 0 & -2^{k} \\ 
-2^{k} & 0 & 1 & 2^{k} \\ 
0 & 2^{k} & 2^{k} & 1%
\end{array}%
\right) $.
\end{quote}

\noindent From the Mayer-Vietoris sequence one gets (e.g. \cite{DC})

\bigskip

\noindent \textbf{Lemma 5.2. }\textsl{Let }$M$\textsl{\ be one of the five
manifolds in (5.1).}

\textsl{i)} \textsl{The second cohomology }$H_{2}(M)$\textsl{\ is}

\begin{quote}
$\qquad H_{2}(S^{2}\times S^{3})=H_{2}(S^{2}\widetilde{\times }S^{3})=%
\mathbb{Z}$; $H_{2}(W)=\mathbb{Z}_{2}$;

$\qquad H_{2}(M_{k})=\mathbb{Z}_{k}\oplus \mathbb{Z}_{k}$; $H_{2}(X_{k})=%
\mathbb{Z}_{2^{k}}\oplus \mathbb{Z}_{2^{k}}$;
\end{quote}

\textsl{ii) }$w_{2}(M)\neq 0$\textsl{\ implies that }$M\in \{S^{2}\widetilde{%
\times }S^{3},W,X_{k}\}.$\hfill $\square $

\bigskip

In terms of the five manifolds constructed in (5.1), the classification of
Smale and Barden on the closed $1$-connected $5$-manifolds \cite{B,S} may be
rephrased as follows.

\bigskip

\noindent \textbf{Lemma 5.3.} \textsl{For any 1-connected 5-manifold }$N\neq
S^{5}$\textsl{\ one has}

\begin{enumerate}
\item[(5.2)] $N=\left\{ 
\begin{tabular}{l}
$\#_{r}S^{2}\times S^{3}\#_{1\leq i\leq t}M_{k_{i}}$ \textsl{if} $%
w_{2}(N)\equiv 0$\textsl{;} \\ 
$\#_{r}S^{2}\times S^{3}\#_{1\leq j\leq t}M_{k_{j}}\#H$ \textsl{if} $%
w_{2}(N)\neq 0$\textsl{,}%
\end{tabular}%
\right. $
\end{enumerate}

\noindent \textsl{where} $H\in \{S^{2}\widetilde{\times }S^{3},W,X_{k}\}$%
\textsl{, and where the notation }$\#_{r}N$ \textsl{means} \textsl{the
connected sum of }$r$\textsl{-copies of} $N$\ \textsl{(and so forth).}

\textsl{In addition, for any primitive element }$\omega \in H^{2}(N)$\textsl{%
\ one can modify the decomposition (5.2) as}

\begin{enumerate}
\item[(5.3)] $N=\left\{ 
\begin{tabular}{l}
$S^{2}\times S^{3}\#N^{\prime }$\textsl{,} \textsl{if} $\omega \neq w_{2}(N)%
\func{mod}2$\textsl{;} \\ 
$S^{2}\widetilde{\times }S^{3}\#N^{\prime }$ \textsl{with} $w_{2}(N^{\prime
})\equiv 0$\textsl{,} \textsl{if} $\omega \equiv w_{2}(N)\func{mod}2$\textsl{%
,}%
\end{tabular}%
\right. $
\end{enumerate}

\noindent \textsl{so that }$\omega $\textsl{\ is a generator of the first
summand of} $H^{2}(N)=H^{2}(S^{2})\oplus H^{2}(N^{\prime })$\textsl{.}\hfill 
$\square $

\bigskip

\noindent \textbf{Proof of Theorem C.} For an $1$-connected $6$-manifold $M$%
\ with a regular circle action let $\pi :$\textsl{\ }$M\rightarrow N:=M/S^{1}
$ be the quotient map. Since the quotient $N$ is then an $1$-connected and $5
$-dimensional manifold which is not diffeomorphic to $S^{5}$, formula (5.2)
is applicable to decompose it as a connected sum of the basic ones.
Moreover, since the Euler class $\omega \in H^{2}(N)$ of $\pi $ is primitive
by Lemma 5.1, we can take a decomposition of $N$ relative to $\omega $ as
that stated in (5.3), to get by the formula (1.4) in Theorem B the
decomposition

\begin{quote}
$M=S^{3}\times S^{3}\#\Sigma _{i}N^{\prime }$, $N^{\prime }\in \mathcal{T}%
_{5}$,
\end{quote}

\noindent where $i=0$ or $1$ in accordance to $\omega \equiv w_{2}(N)\func{%
mod}2$, or $\omega \neq w_{2}(N)\func{mod}2$. Decomposing $N^{\prime }$ in
the form of (5.2) and using (4.1) to expand $\Sigma _{i}N^{\prime }$, we
obtain formula (1.7) in Theorem C, in which 

\begin{quote}
$\Sigma _{0}(S^{2}\times S^{3})=\Sigma _{1}(S^{2}\times S^{3})=S^{3}\times
S^{3}\#S^{2}\times S^{4}$ (by (3.2)).
\end{quote}

\noindent This completes the proof.\hfill $\square $

\section{Torus bundles over 4-manifolds}

The connected-sum operation furnishes the set $\mathcal{T}_{n}$ with the
structure of a semi-abelian group with the sphere $S^{n}$ as the null
element. Write $\mathcal{S}_{n}\subset \mathcal{T}_{n}$ to denote the
semi-subgroup generated by the products $S^{r}\times S^{n-r}$ of spheres.
Since every $M\in $ $\mathcal{S}_{n}$ has the unique decomposition

\begin{quote}
$M=\#_{a_{1}}(S^{1}\times S^{n-1})\#_{a_{2}}{\small \cdots }\#_{a_{[\frac{n}{%
2}]}}(S^{[\frac{n}{2}]}\times S^{n-[\frac{n}{2}]})$,
\end{quote}

\noindent where the sequence of integers $\{a_{1},\cdots ,a_{[\frac{n}{2}]}\}
$ is determined by the Poincar\`{e} polynomial of $M$ as

\begin{quote}
$P_{t}(M)=1+a_{1}(t^{1}+t^{n-1})+\cdots +a_{\left[ \frac{n}{2}\right] }(t^{%
\left[ \frac{n}{2}\right] }+t^{n-\left[ \frac{n}{2}\right] })+t^{n}$,
\end{quote}

\noindent we obtain

\bigskip

\noindent \textbf{Lemma 6.1.} \textsl{The map }$\mathcal{S}_{n}\rightarrow 
\mathbb{Z}[t]$\textsl{\ by }$M\rightarrow P_{t}(M)$ \textsl{injects.}\hfill $%
\square $

\bigskip

\noindent \textbf{Example 6.2.} For $k\geq 1$ the manifold $Q_{k}$
constructed in (1.9) satisfies that

\begin{enumerate}
\item[(6.1)] $P_{t}(Q_{k})=1+c_{1}(t^{3}+t^{k+1})+c_{2}(t^{4}+t^{k})\cdots
+c_{r}(t^{r+2}+t^{k-r-2}),$ $r=\left[ \frac{k}{2}\right] $,
\end{enumerate}

\noindent where the sequence of integers $(c_{1},\cdots ,c_{r})$ is
determined by the integer $k$ by

\begin{quote}
$c_{i}:=(k-1)\binom{k-1}{i}-\binom{k-1}{i+1}+(k-1)\binom{k-1}{i-1}-\binom{k-1%
}{i-2}$, $1\leq i<r$

$c_{r}:=\frac{1}{2}((k-1)\binom{k-1}{r}-\binom{k-1}{r+1}+(k-1)\binom{k-1}{r-1%
}-\binom{k-1}{r-2})$.\hfill $\square $
\end{quote}

The proof of Theorem C (see in Section 1) relies on the following two
lemmas, whose proofs will be postponed.

\bigskip

\noindent \textbf{Lemma 6.3. }\textsl{Let }$M$\textsl{\ be a manifold whose
second cohomology }$H^{2}(M)$\textsl{\ is torsion free with basis }$\left\{
\alpha _{1},\cdots ,\alpha _{k}\right\} $\textsl{. Then, in the notation of
(1,8),}

\textsl{i) If }$\left\{ \alpha _{1}^{\prime },\cdots ,\alpha _{k}^{\prime
}\right\} $ \textsl{is another basis} \textsl{of }$H^{2}(M)$\textsl{, }$%
M(\alpha _{1},\cdots ,\alpha _{k})\cong M(\alpha _{1}^{\prime },\cdots
,\alpha _{k}^{\prime })$\textsl{;}

\textsl{ii) If }$n\geq 5$\textsl{\ and }$M\in \mathcal{S}_{n}$\textsl{,} $%
M(\alpha _{1},\cdots ,\alpha _{k})\in \mathcal{S}_{n+k}$\textsl{.}

\bigskip

\noindent \textbf{Lemma 6.4. }\textsl{If }$N:=\#_{1\leq i\leq
k-1}(S^{2}\times S^{3})$\textsl{\ and if }$\left\{ \omega _{1},\cdots
,\omega _{k-1}\right\} $\textsl{\ is the basis of }$H^{2}(N)$\textsl{\
Kronecker dual to the canonical embeddings }$S^{2}\subset N$\textsl{, then }$%
N(\omega _{1},\cdots ,\omega _{k-1})=Q_{k}\QTR{sl}{.}$

\bigskip

\noindent \textbf{Proof of Theorem D. }Let\textbf{\ }$M$ be an $1$-connected 
$4$-manifold with $\beta _{2}(M)=k$. According to \cite[Corollary 3]{DL}
there is a primitive\textsl{\ }element $e\in H^{2}(M)$ so that the total
space of the circle bundle $\pi :M(e)\rightarrow M$ with Euler class $e$ is

\begin{quote}
$M(e)=\#_{1\leq i\leq k-1}(S^{2}\times S^{3})$ ($=N$ in the notation of
Lemma 6.4).
\end{quote}

\noindent Moreover, in term of the Gysin sequence of $\pi $ \cite[p.143]%
{Mil2}, one finds a basis $\left\{ \alpha _{1},\cdots ,\alpha _{k}\right\} $
of $H^{2}(M)$ such that

\begin{quote}
$\alpha _{1}=e$, $\pi ^{\ast }(\alpha _{i})=\omega _{i-1}$, $2\leq i\leq k$,
\end{quote}

\noindent and that $\{\omega _{1},\cdots ,\omega _{k-1}\}$ is the basis of $%
H^{2}(M(e))$ specified in Lemma 6.4. That is, in the notation of (1.8),

\begin{quote}
$M(\alpha _{1},\cdots ,\alpha _{k})=$ $M(e)(\omega _{1},\cdots ,\omega
_{k-1})=N(\omega _{1},\cdots ,\omega _{k-1})$.
\end{quote}

\noindent We get $M(\alpha _{1},\cdots ,\alpha _{k})=Q_{k}$ by Lemma 6.4.

Finally, if $\left\{ \alpha _{1}^{\prime },\cdots ,\alpha _{k}^{\prime
}\right\} $\textsl{\ }is an arbitrary basis of\textsl{\ }$H^{2}(M)$, we have
by i) of Lemma 6.3 a diffeomorphism 

\begin{quote}
$M(\alpha _{1}^{\prime },\cdots ,\alpha _{k}^{\prime })\cong M(\alpha
_{1},\cdots ,\alpha _{k})$. 
\end{quote}

\noindent This completes the proof of Theorem D.\hfill $\square $

\bigskip

It remains to show Lemmas 6.3 and 6.4.

\bigskip

\noindent \textbf{Proof of Lemma 6.3. }Let $\xi _{k}$ be the universal
principal $T^{k}$-bundle on $BT^{k}$, and let $M$\ be a manifold whose
second cohomology $H^{2}(M)$\ is torsion free. For two bases $\left\{ \alpha
_{1},\cdots ,\alpha _{k}\right\} $ and $\left\{ \alpha _{1}^{\prime },\cdots
,\alpha _{k}^{\prime }\right\} $ of $H^{2}(M)$ with classifying maps $f$ and 
$f^{\prime }:M\rightarrow BT^{k}$, respectively, we have by (1.8) that

\begin{enumerate}
\item[(6.2)] $M(\alpha _{1},\cdots ,\alpha _{k})=f^{\ast }\xi _{k}$, $%
M(\alpha _{1}^{\prime },\cdots ,\alpha _{k}^{\prime })=f^{\prime \ast }\xi
_{k}$.
\end{enumerate}

\noindent On the other hand, expressing the basis elements $\alpha _{i}$'s
in the $\alpha _{i}^{\prime }$'s we get a $k\times k$ integral invertible
matrix $C=(c_{ij})_{k\times k}$ such that

\begin{quote}
$\left\{ \alpha _{1},\cdots ,\alpha _{k}\right\} =\left\{ \alpha
_{1}^{\prime },\cdots ,\alpha _{k}^{\prime }\right\} \cdot C$.
\end{quote}

\noindent Regard $S^{1}$ as the group of the unit complexes $z\in \mathbb{C}$%
, define in term of $C$ the group isomorphism

\begin{quote}
$g:T^{k}\rightarrow T^{k}$, $g(z_{1},\cdots ,z_{k})=(z_{1}^{c_{11}}\cdots
z_{k}^{c_{1k}},\cdots ,z_{1}^{c_{k1}}\cdots z_{k}^{c_{kk}})$,
\end{quote}

\noindent and let $B_{g}:BT^{k}\rightarrow BT^{k}$ be the induced map of $g$%
. From $f^{\ast }=f^{\prime \ast }B_{g}^{\ast }$ we get $f\simeq B_{g}\circ
f^{\prime }$. By the homotopy invariance of the induced bundles, together
with the fact that the map $B_{g}$ has a homotopy inverse, we obtain by
(6.2) a diffeomorphism $M(\alpha _{1},\cdots ,\alpha _{k})\cong M(\alpha
_{1}^{\prime },\cdots ,\alpha _{k}^{\prime })$, showing part i) of Lemma 6.3.

For ii) we can assume by $M\in \mathcal{S}_{n}$ and $n\geq 5$ that

\begin{enumerate}
\item[(6.3)] $M=\#_{k}(S_{i}^{2}\times S_{i}^{n-2})\#_{1\leq j\leq
r}(S_{j}^{t_{j}}\times S_{j}^{n-t_{j}})$, $t_{j}\leq n-t_{j}$, $t_{j}\neq 2$.
\end{enumerate}

\noindent where $S_{t}^{p}\times S_{t}^{n-p}=S^{p}\times S^{n-p}$. For each $%
1\leq i\leq k$ let $q_{i}:M\rightarrow S_{i}^{2}\times S_{i}^{n-2}$ be the
quotient of $M$ by the summands other than the component $S_{i}^{2}\times
S_{i}^{n-2}$, $p_{1}:$ $S_{i}^{2}\times S_{i}^{n-2}\rightarrow S^{2}$ be the
projection onto the first factor, and set

\begin{quote}
$g_{i}=:p_{1}\circ q_{i}:M\rightarrow S_{i}^{2}\times S_{i}^{n-2}\rightarrow
S^{2}$, $\omega _{i}:=g_{i}^{\ast }(\omega )\in H^{2}(M)$,
\end{quote}

\noindent where $\omega \in H^{2}(S^{2})=\mathbb{Z}$ is the orientation
class of $S^{2}$. In term of i) it suffices to show $M(\omega _{1},\cdots
,\omega _{k})\in \mathcal{S}_{n+k}$. This can be done by induction on $k$.

Let $\pi _{1}$ be the circle bundle $M(\omega _{1})\rightarrow M$\ over $M$
with Euler class $\omega _{1}\in H^{2}(M)$. Then

\begin{quote}
$M(\omega _{1},\cdots ,\omega _{k})=M(\omega _{1})(\pi _{1}^{\ast }(\omega
_{2}),\cdots ,\pi _{1}^{\ast }(\omega _{k-1}))$,
\end{quote}

\noindent where, by the Gysin sequence \cite[p.143]{Mil2} of spherical
fibrations, the set $\{\pi _{1}^{\ast }(\omega _{2}),\cdots ,$ $\pi
_{1}^{\ast }(\omega _{k-1})\}$ is a basis of $H^{2}(M(\omega _{1}))$.
However, it follows from Example 1.5, together with formulae (6.3) and
(3.2), that

\begin{center}
$M(\omega _{1})=(S^{3}\times S^{n-2})\#_{k-1}\Sigma _{1}(S_{i}^{2}\times
S_{i}^{n-2})\#_{r}\Sigma _{1}(S_{j}^{k_{j}}\times S_{j}^{n-k_{j}})\in 
\mathcal{S}_{n+1}$.
\end{center}

\noindent This completes the inductive procedure of showing part ii).\hfill $%
\square $

\bigskip

For a finitely generated\textbf{\ }graded\textbf{\ }abelian group $%
A=A^{0}\oplus A^{1}\oplus A^{2}\oplus \cdots $ define the Poincar\`{e}
polynomial $P_{t}(A)$ of $A$ by

\begin{quote}
$P_{t}(A):=a_{0}+a_{1}t+a_{2}t^{2}+\cdots $, $a_{i}=\dim A^{i}\otimes 
\mathbb{R}$.\ \ \ \ \ \ \ \ \ \ \ \ \ \ \ \ \ \ \ \ \ \ \ \ \ \ \ \ \ \ \ \
\ \ \ \ \ \ \ \ \ \ \ \ \ \ \ \ \ \ \ \ \ \ \ \ \ \ \ \ \ \ \ \ \ \ \ \ \ \
\ \ \ \ \ \ \ \ \ \ \ \ \ \ \ \ \ \ \ \ \ \ \ \ \ \ \ \ \ \ \ \ \ \ \ \ \ \
\ \ \ \ \ \ \ \ \ \ \ \ \ \ \ \ \ \ \ \ \ \ \ \ \ \ \ \ \ \ \ \ \ \ \ \ \ \
\ \ \ \ \ \ \ \ \ \ \ \ \ \ \ \ \ \ \ \ \ \ \ \ \ \ \ \ \ \ \ \ \ \ \ \ \ \
\ \ \ \ \ \ \ \ \ \ \ \ \ \ \ \ \ \ \ \ \ \ \ \ \ \ \ \ \ \ \ \ \ \ \ \ \ \
\ \ \ \ \ \ \ \ \ \ \ \ \ \ \ \ \ \ \ \ \ \ \ \ \ \ \ \ \ \ \ \ \ \ \ \ \ \
\ \ \ \ \ \ \ \ \ \ \ \ \ \ \ \ \ \ \ \ \ \ \ \ \ \ \ \ \ \ \ \ \ \ \ \ \ \
\ \ \ \ \ \ \ \ \ \ \ \ \ \ \ \ \ \ \ \ \ \ \ \ \ \ \ \ \ \ \ \ \ \ \ \ \ \
\ \ \ \ \ \ \ \ \ \ \ \ 
\end{quote}

\bigskip

\noindent \textbf{Proof of Lemma 6.4.} Suppose that $N=\#_{1\leq i\leq
k-1}(S_{i}^{2}\times S_{i}^{3})$, and that $\{\omega _{1},\cdots ,\omega
_{k-1}\}$ is the basis of $H^{2}(N)$ specified as that in Lemma 6.4. Then
the cohomology $H^{\ast }(N)$\ has the basis

\begin{enumerate}
\item[(6.4)] $\left\{ 1,\omega _{1},\cdots ,\omega _{k-1},y_{1},\cdots
,y_{k-1},z\right\} $, $\deg \omega _{i}=2$, $\deg y_{i}=3$, $\deg z=5$,
\end{enumerate}

\noindent that is subject to the multiplicative relations

\begin{enumerate}
\item[(6.5)] $\omega _{i}\cdot \omega _{j}=0;$\quad $y_{i}\cdot
y_{j}=0;\quad \omega _{i}\cdot y_{j}=\delta _{ij}\cdot z$, $1\leq i,j\leq
k-1 $,
\end{enumerate}

\noindent where $\delta _{ij}$\ is the Kronecker symbol. By ii) of Lemma
6.3, for the principal $T^{k-1}$-bundle

\begin{quote}
$\pi :N(\omega _{1},\cdots ,\omega _{k-1})\rightarrow N$
\end{quote}

\noindent over $N$ we have $N(\omega _{1},\cdots ,\omega _{k-1})\in S_{4+k}$%
. According to Lemma 6.1 it suffices to derive that

\begin{enumerate}
\item[(6.6)] $P_{t}(N(\omega _{1},\cdots ,\omega _{k-1}))=P_{t}(Q_{k})$ (see
(6.1)).
\end{enumerate}

\noindent To this end we compute with the Leray-Serre spectral sequence $%
\{E_{r}^{\ast ,\ast },d_{r}\}$ of the bundle $\pi $ \cite[p.645]{Wh}.

Since both cohomologies $H^{\ast }(N)$ and $H^{\ast }(T^{k-1})$ are torsion
free, while the base manifold $N$ is $1$-connected, the Leray-Serre theorem
tells that 

\begin{quote}
$E_{2}^{\ast ,\ast }=H^{\ast }(N)\otimes H^{\ast }(T^{k-1})$. 
\end{quote}

\noindent In addition, there are unique elements $t_{1},\cdots ,t_{k-1}\in
H^{1}(T^{k-1})$ such that

\begin{quote}
$H^{\ast }(T^{k-1})=\Lambda (t_{1},\cdots ,t_{k-1})$; $d_{2}(1\otimes
t_{i})=\omega _{i}\otimes 1$,
\end{quote}

\noindent implying

\begin{enumerate}
\item[(6.7)] $d_{2}(x\otimes t_{i})=\omega _{i}\cdot x\otimes 1$, $%
d_{2}(x\otimes 1)=0$, $x\in H^{\ast }(N)$.
\end{enumerate}

For a multi-index $I\subseteq \{1,\cdots ,k-1\}$ we set $t_{I}:=\Pi _{i\in
I}t_{i}$, and let $B_{i}\subset E_{2}^{\ast ,\ast }$, $1\leq i\leq 4$, be
the subgroups with the basis

\begin{quote}
$B_{1}:\{1\otimes t_{I}\mid $ $I\subseteq \{1,\cdots ,k-1\}\}$;$\qquad $

$B_{2}:\{1,\omega _{i}\otimes 1,\omega _{i}\otimes t_{I}\mid \quad i\in
\{1,\cdots ,k-1\},I\subseteq \{1,\cdots ,k-1\}\}$;

$B_{3}:\{y_{i}\otimes 1,y_{i}\otimes t_{I}\mid i\in \{1,\cdots
,k-1\},I\subseteq \{1,\cdots ,k-1\}\}$;$\qquad $

$B_{4}:\{z\otimes 1,z\otimes t_{I}\mid I\subseteq \{1,\cdots ,k-1\}\}$,
\end{quote}

\noindent respectively. Then $E_{2}^{\ast ,\ast }=B_{1}\oplus B_{2}\oplus
B_{3}\oplus B_{4}$ by (6.4) and (as is clear)

\begin{quote}
$P_{t}(B_{1})=\underset{1\leq i\leq k-1}{\Sigma }\binom{k-1}{i}t^{i}$;

$P_{t}(B_{2})=1+(k-1)\underset{0\leq i\leq k-1}{\Sigma }\binom{k-1}{i}%
t^{i+2} $;

$P_{t}(B_{3})=(k-1)\underset{0\leq i\leq k-1}{\Sigma }\binom{k-1}{i}t^{i+3}$;

$P_{t}(B_{4})=\underset{0\leq i\leq k-1}{\Sigma }\binom{k-1}{i}t^{i+5}$.
\end{quote}

\noindent Moreover, the $d_{2}$-action on the basis elements $x\in B_{i}$
has been decided by (6.5) and (6.7) as

\begin{center}
$d_{2}(x)=\left\{ 
\begin{tabular}{l}
$0$ if $x=\omega _{i}\otimes 1,y_{i}\otimes 1,z\otimes t_{I},\omega
_{i}\otimes t_{I},$ $y_{j}\otimes t_{J}$ if $j\notin J$; \\ 
$\underset{1\leq s\leq t}{\Sigma }(-1)^{s-1}\omega _{i_{s}}\otimes t_{I_{s}}$
if $x=1\otimes t_{I}$ with $I=\{i_{1},\cdots ,i_{t}\}$; \\ 
$(-1)^{s-1}z\otimes t_{I_{s}}$ if $x=y_{i}\otimes t_{I}$ with $i=i_{s}\in
I=\{i_{1},\cdots ,i_{t}\}$,%
\end{tabular}%
\right. $
\end{center}

\noindent where $I_{s}$ denotes the multi-index obtained from $I$ by
deleting the $s$-th coordinate $i_{s}$. It follows that

\begin{quote}
a) $d_{2}(B_{2})=0$, $d_{2}$ maps $B_{1}$ monomorphically into $B_{2}$;

b) $d_{2}(B_{4})=0$, $d_{2}(B_{3})\subset B_{4}$ with $B_{4}/d_{2}(B_{3})=%
\mathbb{Z}$ generated by $z\otimes t.$
\end{quote}

\noindent As results,

\begin{enumerate}
\item[(6.8)] $E_{3}^{\ast ,\ast }=B_{2}/d_{2}(B_{1})\oplus
B_{4}/d_{2}(B_{3})\oplus \overline{B}_{3}$ with $\overline{B}_{3}:=\ker
[d_{2}:B_{3}\rightarrow B_{4}]$;

\item[(6.9)] $E_{3}^{\ast ,\ast }=E_{\infty }^{\ast ,\ast }=H^{\ast
}(N(\omega _{1},\cdots ,\omega _{k-1}))$.
\end{enumerate}

\noindent Since

\begin{quote}
$P_{t}(B_{2}/d_{2}(B_{1}))=P_{t}(B_{2})-t\cdot P_{t}(B_{1})$ by a);

$P_{t}(B_{4}/d_{2}(B_{3}))=t^{k+4}$ by b),
\end{quote}

\noindent and since

\begin{quote}
$P_{t}(\overline{B}%
_{3})=P_{t}(B_{3})-t^{-1}(P_{t}(B_{4})-P_{t}(B_{4}/d_{2}(B_{3}))$
\end{quote}

\noindent by the exactness of the sequence of the free $\mathbb{Z}$-modules

\begin{quote}
$\quad 0\rightarrow \overline{B}_{3}\rightarrow B_{3}\overset{d_{2}}{%
\rightarrow }B_{4}\rightarrow B_{4}/d_{2}(B_{3})(=\mathbb{Z)}\rightarrow 0$,
\end{quote}

\noindent we get, in the order of (6.9), (6.8) and (6.1), the equalities

\begin{quote}
$P_{t}(N(\omega _{1},\cdots ,\omega _{k-1}))=P_{t}(E_{3}^{\ast ,\ast })$

$=P_{t}(B_{2})-t\cdot P_{t}(B_{1})+P_{t}(B_{3})-t^{-1}\cdot
(P_{t}(B_{4})-t^{k+4})+t^{k+4}$

$=P_{t}(Q_{k})$.
\end{quote}

\noindent This verifies the relation (6.6), completing the proof of Lemma
6.4.\hfill $\square $

\section{Remarks}

The constructions, results and calculations of the present paper are likely
to be generalized. We discuss such examples.

\textbf{7.1.} By Duan-Liang \cite[Corollary 2]{DL} and Theorem C, the
classification of the $1$-connected $5$- and $6$-dimensional manifolds that
admit regular circle actions has been completed. In addition,
Montgemory-Yang \cite{MY} showed that among all the $28$\ homotopy $7$%
-spheres there are precisely $10$\ admit regular circle actions, which has
been extended by Jiang \cite{J} to the classification of regular circle
actions on the $2$-connected $7$-manifolds, both in the topological and
smooth categories. As suggested by the proofs of Theorem C and Theorem D,
the suspension operations, combined with Theorem B, could provide a powerful
tool to deal with the classification problems in the higher dimensional
cases.

For instance, let $\mathcal{V}_{2n}\subset \mathcal{T}_{2n}$ (resp. $%
\mathcal{V}_{2n+1}\subset \mathcal{T}_{2n+1}$) be the subset of $(n-1)$%
-connected $2n$-manifolds (resp. the subset of $(n-1)$-connected $(2n+1)$%
-manifolds). The operators $\Sigma _{i}$ clearly satisfy the relation $%
\Sigma _{i}(\mathcal{V}_{2n})\subset \mathcal{V}_{2n+1}$. On the other hand,
the manifolds of $\mathcal{V}_{m}$ ($m\geq 8$) have been classified by
C.T.C. Wall \cite{Wall1,Wall2} in term of a set of explicit invariants. It
is possible to describe the map $\Sigma _{i}:\mathcal{V}_{2n}\rightarrow 
\mathcal{V}_{2n+1}$ in terms of the Wall's invariants, so that a
classification for the manifolds of $\mathcal{V}_{2n+1}$ that admit regular
circle actions is attainable.

\textbf{7.2.} The operations $\Sigma _{i}:\mathcal{T}_{n}\rightarrow 
\mathcal{T}_{n+1}$ admit a natural generalization. Namely, for each pair $%
(\alpha ,M)\in \pi _{k}(SO(n))\times \mathcal{T}_{n}$ take an orientation
persevering embedding $(D^{n},0)\subset (M,x_{0})$, surgery along the
oriented sphere $S^{k}\times \{x_{0}\}\subset S^{k}\times M$ relative to the
framing $\alpha \in \pi _{k}(SO(n))$, and denote the resulting manifold by $%
\Sigma _{\alpha }M$ to get the pairing

\begin{quote}
$\Sigma :\pi _{k}(SO(n))\times \mathcal{T}_{n}\rightarrow \mathcal{T}_{n+k}$%
, $(\alpha ,M)\rightarrow \Sigma _{\alpha }M$.
\end{quote}

\noindent As indicated by the construction (1.3), this operation may be use
to describe the total space $E_{N}$ the $S^{k}$-bundles $p_{N}:E_{N}%
\rightarrow B\#N$ induced from a $S^{k}$-bundles $E\rightarrow B$, where $%
\dim B=\dim N$. In particular, when $k=3$ it may be applied to construct or
classify the $S^{3}=SU(2)$-actions on smooth manifolds.

\bigskip

Haibao Duan

dhb@math.ac.cn

Yau Mathematical Science Center, Tsinghua University, Beijing 100084;

Academy of Mathematics and Systems Sciences, Chinese Academy of Sciences,
Beijing 100190.

\end{document}